\documentclass[11pt]{amsart}
\usepackage[T1]{fontenc}
\usepackage{a4wide}
\usepackage{psfrag}
\usepackage[dvips]{graphicx}
\usepackage[english]{babel}
\usepackage{amsmath}
\usepackage{amsthm}
\usepackage{amssymb}
\usepackage{amsfonts}
\usepackage{fullpage}
\usepackage{mathrsfs}
\usepackage[all]{xy}

\newcommand\C{\mathbb C}

\newcommand\Q{\mathcal Q}

\newcommand\CP{{\mathbb C}\text{P}^1}
\newcommand\D{d\hspace{-0.5pt}}

\renewcommand{\mod}{\:\operatorname{mod}}

\newtheorem{Theorem}{Theorem}
\newtheorem{Corollary}{Corollary}
\newtheorem{Proposition}{Proposition}
\newtheorem{Lemma}{Lemma}
\newtheorem{Construction}{Construction}

\newtheorem*{NoNumberTheorem}{Theorem}
\newtheorem*{NoNumberProposition}{Proposition}
\newtheorem*{Riemann:Hurwitz}{Riemann--Hurwitz Formula}

\theoremstyle{remark}
\newtheorem{Remark}{Remark}

\theoremstyle{definition}
\newtheorem{Definition}{Definition}

\newtheorem*{Convention}{Convention}
\newtheorem*{Forgetful Map}{``Forgetful'' Map}
\newtheorem*{Cohomological coordinates}{Cohomological coordinates}
\newtheorem*{Caseone}{The case $n+r = 0$ and $m=0$}
\newtheorem*{Casetwo}{The case $n+r=1$ and $m \geq 0$}
\newtheorem*{Casethree}{The case $n+r=0$ and $m > 0$}

\begin{document}
\title[Hyperelliptic components of the moduli spaces]
{Hyperelliptic Components of the Moduli Spaces of Quadratic
Differentials with Prescribed Singularities}

\author{Erwan Lanneau}

\address{Institut de math\'ematiques de Luminy (IML)\\
163 Avenue de Luminy, Case 907 \\
F-13288 Marseille Cedex 9, France}

\email{lanneau@iml.univ-mrs.fr}

\subjclass{32G15, 30F30, 30F60, 58F18}
\keywords{Quadratic differentials, Moduli space, Measured foliations, 
Teich\"muller geodesic flow}
\date{\today}

\begin{abstract}
Moduli spaces of quadratic differentials with prescribed
singularities are not necessarily connected. We describe here all
cases when they have a special {\it hyperelliptic} connected component.

We announce the general classification theorem: up to the four 
exceptional cases in low dimensional stratum, any stratum of meromorphic 
quadratic differentials is either connected, or has exactly two 
connected components. In this last case, one component is hyperelliptic, 
the other not.
\end{abstract}

\maketitle

\setcounter{tocdepth}{1}
\tableofcontents

\section {Introduction}

A meromorphic quadratic differential $\psi$ on a Riemann
surface $M^2_g$ of genus $g$ is locally defined by $\psi=f(z)(\D
z)^2$ where $f(z)$ is a meromorphic function defined on a chart $(U,z)$. 
If we have the form $\psi=g(w)(\D w)^2$ where $g(w)$ is a 
meromorphic function on $V$, the transition functions 
on $U \cap V$ satisfies
$$
\cfrac{f(z)}{g(w)}=\left( \cfrac{\D w}{\D z} \right)^2
$$
In this paper we
consider quadratic differentials having only simple poles, if any.
We denote by $\Q_g$ the moduli spaces of pairs $(M^2_g,\psi)$,
where $M^2_g$ is a Riemann surface and $\psi$ a meromorphic
quadratic differential on it. Note that we have divided by 
the mapping class group Mod$(g)$.

Let $(k_1,\dots,k_n)$ be the orders of singularities of $\psi$,
where $k_i>0$ corresponds to a zero $z^{k_i}(\D z)^2$ of order
$k_i$ and $k_j=-1$ corresponds to a simple pole $\cfrac{1}{z}\,
(\D z)^2$ of a quadratic differential. It is a classical result 
following of the Gauss---Bonnet formula that $\sum k_i = 4g-4$.

Locally, in a simply-connected neighborhood of a nonsingular point, 
a quadratic differential can be presented as a square of an
Abelian differential, but globally it is not the case in general.
In this paper we consider only those quadratic differentials which
are not globally the squares of Abelian differentials.

The moduli space $\Q_g$ is naturally stratified by the types of the
singularities. We denote by $\Q(k_1,\dots,k_n)\subseteq \Q_g$ the
stratum of quadratic differentials $[M^2_g,\psi]\in\Q_g$ which are
not the global squares of Abelian differentials,  and which have the
singularity pattern $(k_1,\dots,k_n)$, where $k_i$ takes values in 
$\{-1,0,1,2,\dots\}$. In some situations we consider the strata of
pairs $[\psi, M^2_g]$ where the Riemann surface is provided with
some fixed number of marked points, which are nonsingular points
of the quadratic differential $\psi$. By convention we let $k_i=0$
when the point $P_i$ is, actually, just a marked point. Sometimes we
shall call the marked points the {\it fake} zeros, and the zeros
$P_j$ with $k_j>0$ the {\it true} zeros. We use the exponential
notation $k^m$ for $k,k,k,\dots,k$ repeated $m$ times. For
example, $\Q(1^4,8,2,3^2)$ stands for $\Q(1,1,1,1,8,2,3,3)$.

We are interested in the description of connected components 
of the moduli spaces $\Q(k_1,\dots,k_n)$. In this paper we describe 
all connected components of a special type: the {\it hyperelliptic}
ones. A complete description of all connected components of the
moduli spaces $\Q(k_1,\dots,k_n)$ can be found 
in~\cite{Lanneau:02:classification}. It will be the
subject of a forthcoming paper.

\begin{Remark}
In~\cite{Kontsevich:Zorich} Kontsevich and Zorich classified
connected components of the moduli spaces of Abelian
differentials, which implies the classification of the moduli
spaces of those holomorphic quadratic differentials which are
globally the squares of Abelian differentials.
\end{Remark}

\begin{Remark}
In~\cite{Masur:82} and~\cite{Veech:82}, Masur and Veech have independently proved
that the Teichm\"uller geodesic flow acts ergodically on each 
connected component of each stratum of the moduli space of quadratic differentials; 
the corresponding invariant measure being a finite Lebesque equivalent measure.

Motivated by this result, the classification of connected components of the 
strata coincide with the classification of ergodic components of the 
Teichm\"uller geodesic flow.
\end{Remark}

The paper has the following structure. In
section~\ref{s:covering:construction} we present an overview of
the basic properties of quadratic differentials. In
section~\ref{s:hyperelliptic:components} we construct the hyperelliptic
connected components and prove Theorem~\ref{thm:main}, that is 
our list is complete. In section
~\ref{s:qd:versa:foliations} we present the well-known relation
between quadratic differentials and measured foliations. Then we use measured
foliations to prove Theorem~\ref{th:non:connected}: 
strata with a hyperelliptic component are not
connected except some particular cases in low genera. We draw from our 
two main results the following

\begin{Corollary}
\textrm{ }

\begin{itemize}

\item The three following series of strata in $\Q_g$
$$
\begin{array}{llll}

\mathcal F_2 = \{ & \Q(4(g-k)-6\ ;\ 4k+2)  & | & \ \ 0 \leq k \leq g-2 \}  \\

\mathcal F_3 = \{ & \Q(4(g-k)-6\ ;\ 2k+1\ ;\ 2k+1) & | & \ \ 0 \leq k \leq g-1 \}  \\

\mathcal F_4 = \{ & \Q((2(g-k)-3\ ;\ 2(g-k)-3\ ;\ 2k+1\ ; 2k+1) &
|  & -1 \leq k  \leq g-2 \}

\end{array}
$$
in genera $g \geq 3$ are non-connected. Each above stratum 
possesses one component which is hyperelliptic and at least one other 
which is not.

\vskip 5mm

\item The following strata corresponding respectively to genus $1$ and $2$
are connected and the whole stratum coincides with its hyperelliptic connected component
\begin{displaymath}
\begin{array}{ccc}
\begin{cases}
\Q(-1,-1,2) \\
\Q(-1,-1,1,1)
\end{cases}  &
and &
\begin{cases}
\Q(2,2) \\
\Q(1,1,2) \\
\Q(1,1,1,1)
\end{cases}
\end{array}
\end{displaymath}
\end{itemize}
\end{Corollary}

In section~\ref{general:theorem} we announce the general classification 
theorem: up to four exceptional cases in low genera the strata 
of meromorphic quadratic differentials are either connected,
or have exactly two connected components,
and one of the two components is hyperelliptic.

\section{Mappings of the Moduli Spaces Induced by Ramified
Coverings of a Fixed Combinatorial Type}
\label{s:covering:construction}

In this section we present some general information concerning the
moduli spaces of quadratic differentials. All proofs and 
details can be found in papers~\cite{Masur:82}, \cite{Veech:82},
\cite{Masur:Smillie}, \cite{Veech}, \cite{Kontsevich},
\cite{Kontsevich:Zorich}.

\begin{NoNumberTheorem}[Masur and Smillie]
Consider a vector $(k_1,\dots,k_n)$ with all $k_i \in \mathbb{N}
\cup \{-1 \}$. Suppose that $\sum k_i = 0 \mod 4$ and $\sum k_i
\geq -4$. Then the corresponding stratum $\Q(k_1,\dots,k_n)$ is
non-empty with the following four exceptions
$$
 \Q(\emptyset), \Q(1,-1)\ (\text{in genus } g=1)
 \quad \text{and}\quad
 \Q(4), \Q(1,3)\ (\text{in genus } g=2)
$$
\end{NoNumberTheorem}

\begin{Remark}
It is clear that if a quadratic differential has a singularity of
odd order ($-1$ for a pole) then it may not be a global square of
an Abelian differential. But if all zeros are of even orders
$k_i$ (except $(k_1,\dots,k_n)=(\emptyset)$ or $(4)$) then the
Theorem above says that there exist quadratic differentials with
the singularity pattern $(k_1,\dots,k_n)$ which are not squares of
Abelian differentials. For example, on each complex curve of genus
$g=2$, there exist quadratic differentials with two zeros of
order two which {\it are not} squares of Abelian differentials
(see also Figure~\ref{fig:one:cylinder}), and there exist
quadratic differentials with two zeros of order two which {\it
are} the global squares of Abelian differentials.
\end{Remark}

\begin{NoNumberTheorem}[H.~Masur; W.~Veech]
Any stratum $Q(k_1,\dots,k_n)$ is a complex orbifold of dimension
\medskip

\centerline {${\rm dim}_{\mathbb{C}} \Q(k_1,\dots,k_n) = 2g + n
-2$}
\end{NoNumberTheorem}

\begin{Remark}
The principal stratum, that is all singularities of the form are only 
simples zeroes, 
can be identified with the cotangent fiber over the Teichm\"uller 
space. So the dimension in this case is 
$$
2 \cdot \textrm{ dim}_\mathbb{C} T_g=6g-6
$$
This correspond to the above formula with $n=4g-4$. The general 
formula is obtained by subtracting a dimension any time 
a zero is collapsed to higher order.
\end{Remark}

\begin{NoNumberProposition}[Kontsevich]
Any stratum $\Q(k_1,\dots,k_n)$ with $\sum k_i = -4$ is connected.
\end{NoNumberProposition}

\begin{proof} Since there is only one complex structure on
$\CP$ we can work in the standard atlas on
$\CP=\C \cup (\C^\ast \cup \infty)$. In
this atlas, we can easily find quadratic differentials $f(z)(\D
z)^2$ with any prescribed singularities at any prescribed points
(with the evident condition prescribed by degrees of singularities, 
$\sum k_i = -4$) just
by choosing an appropriate rational function $f(z)$. The space of
configurations of points on a sphere is connected; this implies
the statement of the Proposition.
\end{proof}

\begin{Construction}[Canonical double covering]
\label{constr:canonical:2:covering}
Let $M^2_g$ be a Riemann surface and let $\psi$ be a quadratic
differential on it which is not a square of an Abelian
differential. There exists a canonical (ramified) double covering
$\pi : \tilde{M}^2_{\tilde{g}} \rightarrow M^2_g$ such that
$\pi^{*} \psi=\tilde{\omega}^2$, where $\tilde{\omega}$ is an
Abelian differential on $\tilde{M}^2_{\tilde{g}}$.

The set of critical values of $\pi$ on $M^2_g$ coincide exactly 
with the set of singularities of odd degrees of $\psi$. The
covering $\pi: \tilde{M}^2_{\tilde{g}} \rightarrow M^2_g$ is the
minimal (ramified) covering such that the quadratic differential
$\pi^\ast\psi$ becomes the square of an Abelian differential on
$\tilde{M}^2_{\tilde{g}}$.
\end{Construction}

\begin{proof} Consider an atlas
$(U_i,z_i)_i$ on $\dot{M}^2_g=M^2_g \backslash
\{\text{singularities of }\psi\}$ where we punctured all the zeros
and the poles of $\psi$. We assume that all the charts $U_i$ are
connected and simply-connected. The quadratic differential $\psi$
can be represented in this atlas by a collection of holomorphic
functions $f_i(z_i)$, where $z_i\in U_i$, satisfying the
relations:
$$
 f_i\left(z_i(z_j)\right)\cdot\left(\cfrac{\D z_i}{\D z_j}\right)^2 =
 f_j(z_j)\ \text{ on }U_i \cap U_j
$$
Since we have punctured all singularities of $\psi$ any function
$f_i(z_i)$ is nonzero at $U_i$. Consider two copies $U^{\pm}_i$ of
every chart $U_i$: one copy for every of two branches
$g_i^\pm(z_i)$ of $g^\pm(z_i):=\sqrt{f_i(z_i)}$ (of course, the
assignment of ``$+$'' or ``$-$'' is not canonical). Now  for every
$i$ identify the part of $U^{+}_i$ corresponding to $U_i\cap U_j$
with the part of one of $U^{\pm}_j$ corresponding to $U_j\cap U_i$
in such way that on the overlap branches would match
$$
 g^{+}_i(z_i(z_j))\cdot\cfrac{\D z_i}{\D z_j} =
 g^{\pm}_j(z_j)\ \text{ on }U^{+}_i \cap U^{\pm}_j
$$
Apply the analogous identification to every $U^{-}_i$. We get a
Riemann surface with punctures provided with a holomorphic 1-form
$\tilde{\omega}$ on it, where $\tilde{\omega}$ is presented by the
collection of holomorphic functions $g^{\pm}_i$ in the local
charts. It is an easy exercise to check that filling punctures
we get a closed Riemann surface $\tilde{M}^2_{\tilde{g}}$, and
that $\tilde{\omega}$ extends to an Abelian differential on it. We
get a canonical (possibly ramified) double
covering $\pi:\tilde{M}^2_{\tilde{g}}\to M^2_g$ such that
$\pi^{\ast}\psi=\tilde{\omega}^2$.

By construction the only points of the base $M^2_g$ where the
covering might be ramified are the singularities of $\psi$. In a
small neighborhood of a zero of even degree $2k$ of $\psi$ we can
chose coordinates in which $\psi$ is presented as $z^{2k} (\D
z)^2$. In this chart we get two distinct branches $\pm z^k \D z$
of the square root. Thus the zeros of even degrees of $\psi$ and
the marked points are the regular points of the covering $\pi$.
However, it easy to see that the covering $\pi$ has  a
ramification point over any zero of odd degree and over any simple
pole of $\psi$.
\end{proof}

\begin{Cohomological coordinates}
We use the construction above to
describe the local cohomological coordinates on
$\Q(k_1,\dots,k_n)$ proposed by Kontsevich.

Let $[M,\psi]\in \Q(k_1,\dots,k_n)$. Consider the canonical double
covering $\pi:\tilde{M}\to M$ described in above 
Construction~\ref{constr:canonical:2:covering} such that the 
pull-back $\pi^\ast\psi=\omega^2$ becomes the global
square of an Abelian differential $\omega$ on $\tilde{M}$. Let
$\tau$ be the natural involution of $\tilde{M}$ interchanging the
points in the fibers of $\pi$.  Let
$\tilde{P}_1,\dots,\tilde{P}_r\in\tilde{M}$ be the {\it true}
zeros of $\omega$. Since by construction
$\tau^\ast\omega=-\omega$, the set
$\{\tilde{P}_1,\dots,\tilde{P}_r\}$ is sent to itself by the
involution $\tau$. Consider the induced involution
$$
 \tau^\ast: H^1(\tilde{M},\{ \tilde{P}_1, \dots, \tilde{P}_r \};\C) \to
            H^1(\tilde{M},\{ \tilde{P}_1, \dots, \tilde{P}_r \};\C)
 $$
of the relative cohomology group. The vector space
$H^1(\tilde{M},\{ \tilde{P}_1, \dots, \tilde{P}_r \};\C)$ splits
into direct sum
$$
 H^1(\tilde{M},\{ \tilde{P}_1, \dots, \tilde{P}_r \};\C) = V_1 \oplus V_{-1}
$$
of invariant and anti-invariant subspaces of the involution
$\tau^\ast$. We have already seen that $[\omega]\in V_{-1}$. A
small neighborhood of $[\omega]$ in the anti-invariant
subspace $V_{-1}$ gives a local coordinate chart of a regular
point $[M,\psi]$ (not a fixed point for a elliptic element of
Mod$(g)$) in the stratum $\Q(k_1,\dots,k_n)$.\qed
\end{Cohomological coordinates}

In what follows we consider ramified coverings
$\pi:\tilde{M}^2_{\tilde{g}} \to M^2_g $ of arbitrary degree $d$.
We denote the ramification index of $\pi$ at a point
$\tilde{P}\in\tilde{M}^2_{\tilde{g}}$ by $e_\pi(\tilde{P})$. By
convention we have $e_\pi(\tilde{P})=1$ at any regular point
$\tilde{P}$ of the covering $\pi$. The Riemann--Hurwitz formula
gives the value of the genus $\tilde{g}$ of the covering Riemann
surface $\tilde{M}^2_{\tilde{g}}$.

\begin{Riemann:Hurwitz}
Let $\pi:\tilde{M}^2_{\tilde{g}} \to M^2_g $ be an analytic map of
degree $d$ between compact Riemann surfaces. The genus $\tilde{g}$
of $\tilde{M}^2_{\tilde{g}}$ and the genus $g$ of $M^2_g$ are
related by the formula
$$
2\tilde{g}-2=
d\cdot (2g-2)+\sum_{\tilde{P}\in\tilde{M}^2_{\tilde{g}}}(e_\pi(\tilde{P})-1)
$$
where $e_\pi(\tilde{P})$ is the index of the ramification of $\pi$
at $\tilde{P}$.
\end{Riemann:Hurwitz}

Having a ramified $d$-fold covering $\pi:\tilde{M}^2_{\tilde{g}}
\to M^2_g $ and a meromorphic quadratic differential $\psi$ on
$M^2_g$ with a singularity pattern $(k_1,\dots,k_n)$  we will
need to compute the singularity pattern
$(\tilde{k}_1,\dots,\tilde{k}_m)$ of the induced quadratic
differential $\pi^\ast\psi$ on $\tilde{M}^2_{\tilde{g}}$.

\begin{Lemma}
\label{lm:change:of:degree}
Let $\pi:\tilde{M}^2_{\tilde{g}} \to M^2_g $ be a (ramified)
covering, and let $\psi$ be a meromorphic quadratic differential
on $M^2_g$. A point $\tilde{P}\in\tilde{M}^2_{\tilde{g}}$ is a
singular point of the induced quadratic differential
$\pi^\ast\psi$ either if its image $P=\pi(\tilde{P})$ is a
singular point of $\psi$ or if $\tilde{P}$ is a ramification point
of $\pi$. The degree $\tilde{k}$ of $\pi^\ast\psi$ at the point
$\tilde{P}$ and the degree $k$ of $\psi$ at the point
$P=\pi(\tilde{P})$ are related as
$$
\tilde{k}=e_{\pi}(\tilde{P}) \cdot (k+2)-2,
$$
\end{Lemma}
\begin{proof}
Let $P \in M^2_g$ be a point of the base of the covering.
Let the quadratic differential have a singularity of order $k$ at $P$.
Recall that $k=-1$ if $P$ is a simple pole of $\psi$, and $k=0$ if
$P$ is, actually, a regular point. Let $\tilde{P}$ be in the
preimage $\pi^{-1}(P)$ of $P$. If $\tilde{P}$ is a regular point
of the covering $\pi$ then $\pi^{\ast}\psi$ has at $\tilde{P}$ a
singularity of the same degree $k$ as the singularity $P$ of
$\psi$.

Suppose that $\tilde{P}$ is a ramification point of $\pi$ of index
$e_{\pi}(\tilde{P})=b$. We can choose a local coordinate $z$ in
the neighborhood of the point $P\in M^2_g$ in such way that $\psi=
z^k (\D z)^2$ in this coordinate. We can now chose the local
coordinate $w$ in a neighborhood of $\tilde{P}\in
\tilde{M}^2_{\tilde{g}}$ in such way that the projection $\pi$ has
the form $z = w^b$ in this coordinate. Then we get the following
representation of the induced quadratic differential
$\pi^{\ast}\psi$ in the neighborhood of the point $\tilde{P}$
$$
\pi^{\ast}\psi = (z(w))^k \Big( \D (z(w)) \Big)^2 \sim w ^{b\cdot k}
( w^{b-1} \D w)^2 = w^{b\cdot k+2b-2} (\D w)^2
$$

In particular, if the order of the covering is $2$ then a
singularity of $\psi$ gives either one singularity of order $2k+2$
or two singularities of order $k$ of $\pi^\ast \psi$.
\end{proof}

We use the following obvious Corollary of
Lemma~\ref{lm:change:of:degree}.

\begin{Corollary}
\label{cr:disappearing:singularity}
Let $\pi:\tilde{M}^2_{\tilde{g}} \to M^2_g $ be a ramified
covering, and let $\psi$ be a meromorphic quadratic differential
on $M^2_g$. A preimage $\tilde{P}=\pi^{-1}(P)$ of a singular point
$P$ of $\psi$ is a regular point of $\pi^\ast\psi$ only in the
following two cases
\begin{itemize}
\item the ramification index $e_{\pi}(\tilde{P})$ of the point $\tilde{P}$
is equal to $2$, and $\psi$ has a simple pole at $P$;
\item the ramification index $e_{\pi}(\tilde{P})$ of the point $\tilde{P}$
is equal to $1$ and $P$ is a fake zero of $\psi$.
\end{itemize}
\end{Corollary}

\rightline{$\Box$}

To complete this section we present a construction of a natural
mapping of the strata induced by a ramified covering of the fixed
combinatorial type. This mapping was introduced
in~\cite{Kontsevich:Zorich} to construct the hyperelliptic connected
components of the moduli spaces of Abelian differentials. In the
next section we shall use this mapping to construct the hyperelliptic
connected components of the moduli spaces of quadratic
differentials.

Let $M^2_g$ be a Riemann surface and let $\psi_0$ be a quadratic
differential on it which is not a square of an Abelian
differential. Let $(k_1, \dots, k_n)$ be its singularity pattern.
We do not exclude the case when some of $k_i$ are equal to zero:
by convention this means that we have some marked points.

Let $\pi : \tilde{M}^2_{\tilde{g}} \rightarrow M^2_g$ be  a
(ramified) covering such that the image of any ramification point
of $\pi$ is a marked point, or a zero, or a pole of the quadratic
differential $\psi_0$. Fix the {\it combinatorial type} of the
covering $\pi$: the degree of the covering, the number of critical
fibers and the ramification index of the points in every critical
fiber. Consider 
the induced quadratic differential $\pi^{\ast}\psi_0$ on
$\tilde{M}^2_{\tilde{g}}$; let $(\tilde{k}_1,\dots,\tilde{k}_m)$
be its singularity pattern.

\begin{Construction}[Ramified covering construction]
\label{constr:general}
Deforming slightly the initial point
$[M^2_g,\psi_0]\in\Q(k_1,\dots,k_n)$ we can consider a ramified
covering over the deformed Riemann surface of the same
combinatorial type as the covering $\pi$. This new covering has
exactly the same relation between the position and types of
the ramification points and the degrees and position of the 
singularities of the deformed quadratic differential. This means
that the induced quadratic differential $\pi^{\ast}\psi$ 
has the same singularity pattern
$(\tilde{k}_1,\dots,\tilde{k}_m)$ as $\pi^{\ast}\psi_0$. Thus we
get a local mapping
\begin{gather*}
\Q(k_1, \dots, k_n) \to \Q(\tilde{k}_1,\dots,\tilde{k}_m)\\
[M^2_g,\psi]\mapsto (\tilde{M}^2_g,\pi^{\ast}\psi)
\end{gather*}
\end{Construction}

Note that in general the corresponding {\it global} mapping is
multi-valued. For example, if there is a pair of singularities of
$\psi$ of the same order $k$ such that the first one is the image
of a ramification point of the covering and the other one is not,
then the corresponding map is multi-valued.

Note also that for ramified coverings of some
special types the image of the mapping belongs to a stratum of
squares of Abelian differentials. This is the case, for example,
for the canonical double covering described in
Construction~\ref{constr:canonical:2:covering}.

\begin{Forgetful Map}
Consider a stratum of meromorphic quadratic differentials with $m$ marked points
$$\Q(k_1,\dots,k_n,k_{n+1},\dots,k_{n+m})$$ where $k_i \not = 0$ for 
$i=1,\dots,n$ and $k_{i+n}=0$ for $i=1,\dots,m$. 
We have the following forgetful map
$$
\Q(k_1,\dots,k_n,k_{n+1},\dots,k_{n+m}) \to \Q(k_1,\dots,k_n)
$$
\end{Forgetful Map}

\begin{Remark}
We are interested by the topology of the strata, that is the 
classification of connected components of strata $\Q(k_1,\dots,k_n)$ 
with all $k_i \not = 0$. One can prove that the topology of 
strata with marked points and the topology of strata without 
marked point coincide. Nevertheless we use strata with marked point, 
that is some $k_i$ is equal to zero, for the above construction.

We use Construction~\ref{constr:general} to obtain connected components 
of the strata with {\it no} marked points. Thus we use, actually, the 
composition of the mapping defined above with the forgetful map.
\end{Remark}

\begin{Convention}
By convention we consider only those coverings $\pi : \tilde{M} \to M$, 
that have at least one ramification point in the fiber over any marked 
point of $\psi$ on the underlying surface.

-- We shall call the images of ramification points the {\it critical
values} of $\pi$. Our convention means that the marked points of $\psi$ form a
subset of the critical values of $\pi$.
\end{Convention}

We complete this section with the following statement.
\begin{Lemma}
\label{lm:mapping:is:not:degenerate}
With the above convention the mapping
 $$
  \Q(k_1, \dots, k_n) \to \Q(\tilde{k}_1,\dots,\tilde{k}_m)
 $$
is locally an embedding.
\end{Lemma}

\begin{proof}[Proof of Lemma~\ref{lm:mapping:is:not:degenerate}]
We need to prove that the mapping constructed above is locally injective
near a regular point $[M,\psi]$ in $\Q(k_1,\dots,k_n)$. We use the {\it 
cohomological coordinates} (see above)
in the neighborhood of $[M,\psi]\in \Q(k_1,\dots,k_n)$ and in the
neighborhood of its image $[\tilde{M},\tilde{\psi}]\in
\Q(\tilde{k}_1,\dots,\tilde{k}_m)$,
where $\tilde{\psi}=\pi^\ast\psi$. These coordinates
linearize the mapping, and the proof becomes an exercise
in algebraic topology.

Take the two {\it canonical double coverings} (see
Construction~\ref{constr:canonical:2:covering})
$$
p : N \to M \quad\textrm{ and }\quad
\tilde{p} : \tilde{N} \to \tilde{M}
$$
such that the pull-back $p^\ast\psi=\omega^2$ of the quadratic
differential $\psi$ from $M$ to $N$ becomes the global square of
an Abelian differential $\omega$ on $N$, and the pull-back
$\tilde{p}^\ast\tilde{\psi}=\tilde{\omega}^2$ of the quadratic
differential $\tilde{\psi}$ from $\tilde{M}$ to $\tilde{N}$
becomes the global square of an Abelian differential
$\tilde{\omega}$ on $\tilde{N}$. If the quadratic differential
$\tilde{\psi}$ is already the global square of an Abelian
differential on $\tilde{M}$, the double covering $\tilde{p} :
\tilde{N} \to \tilde{M}$ is not connected: it is composed of two
copies of $\tilde{M}$, and the corresponding Abelian differential
on $\tilde{N}$ corresponds to two branches $\pm\,\tilde{\omega}$
of the globally defined square root of $\tilde{\psi}$.

It follows from the definition of the canonical double covering,
that the diagram
$$
\xymatrix{
\tilde{N} \ar[r]^{\tilde{p}} & \tilde{M} \ar[d]^{\pi} \\
N \ar[r]^{p} & M
}
$$
can be completed to a commutative diagram
$$
\xymatrix{
\tilde{N} \ar[d]^{\tilde{\pi}} \ar[r]^{\tilde{p}} & \tilde{M} \ar[d]^{\pi} \\
N \ar[r]^{p} & M
}
$$
Moreover, the induced (ramified) covering
$\tilde{\pi}:\tilde{N}\to N$ can be chosen in such way that
$\tilde{\omega}=\tilde{\pi}^\ast\omega$. In particular, it
intertwines the natural involutions
$\tilde{\tau}:\tilde{N}\to\tilde{N}$ and $\tau: N\to N$, that is
$\tilde{\pi}\circ\tilde{\tau}=\tau\circ\tilde{\pi}$.

Denote by $Q_1, \dots, Q_r$ the zeroes of $\omega$ on
$N$. The set $\{Q_1, \dots, Q_r\}$ is obtained as an inverse image $p^{-1}$
of the set of true zeros, the marked points and the poles of $\psi$.

Note that we do not mark any points at $\tilde{M}$ (we are 
interested by classification of strata without marked point, that 
is without fake singularities). Then the Convention draw aside the 
trivial cases when the map is degenerated.
 
Denote by $\tilde{Q}_1, \dots, \tilde{Q}_l$ the zeros of 
$\tilde{\omega}$ on $\tilde{N}$. 
Since $\tilde{\omega}=\tilde{\pi}^\ast\omega$, we 
see that any zero $\tilde{Q}_i$ of $\tilde{\omega}$ located at a
regular point of the covering $\tilde{\pi}$ is projected by
$\tilde{\pi}$ to a zero of $\omega$. It follows from
Convention that the image of a zero
$\tilde{Q}_j$ of $\tilde{\omega}$ located at a ramification point
of the covering $\tilde{\pi}$ is also projected by $\tilde{\pi}$ to a
zero or to a marked point of $\omega$. Thus
$$
\tilde{\pi}(\{\tilde{Q}_1, \dots, \tilde{Q}_l\})\subset
\{Q_1,\dots, Q_r\}
$$
Hence the covering map $\tilde{\pi}$ induces a mapping
\begin{equation}
\label{eq:tilde:pi:ast}
\tilde{\pi}^\ast : H^1(N,\{ Q_1,\dots,Q_r \},\C) \to
H^1(\tilde{N},\{ \tilde{Q}_1,\dots,\tilde{Q}_l \},\C)
\end{equation}
Since this mapping intertwines the natural involutions
$\tilde{\pi}^\ast\circ\tilde{\tau}^\ast=
\tau^\ast\circ\tilde{\pi}^\ast$, the subset
$V_{-1}\subset H^1(N,\{ Q_1,\dots,Q_r \},\C)$ anti-invariant
under the involution $\tau^\ast$ is mapped into the subset
$\tilde{V}_{-1}\subset
H^1(\tilde{N},\{ \tilde{Q}_1,\dots,\tilde{Q}_l \},\C)$
anti-invariant under the involution $\tilde{\tau}^\ast$.
The induced map $\tilde{\pi}^\ast : V_{-1}\to \tilde{V}_{-1}$ restricted
to an appropriate neighborhoods of $[M,\psi]\in V_{-1}$ and
$[\tilde{M},\tilde{\psi}]\in \tilde{V}_{-1}$ coincides with the
mapping $\Q(k_1,\dots,k_n)\to\tilde{\Q}(\tilde{k}_1,\dots,\tilde{k}_m)$
written down in the cohomological coordinates.

It is obvious that the map of the {\it absolute} cohomology groups
$$
\tilde{\pi}^\ast : H^1(N;\C) \to H^1(\tilde{N};\C)
$$
induced by the covering $\pi : \tilde{N} \to N$ is a monomorphism.
It is not difficult to see that, with our convention on the covering $\pi$,
the restriction of the mapping $\tilde{\pi}^\ast : V_{-1}\to 
\tilde{V}_{-1}$ to the anti-invariant subspace of $\tau^\ast$ is also 
a monomorphism, which completes the proof of Lemma~\ref{lm:mapping:is:not:degenerate}.
Note that in general, the mapping $\tilde{\pi}^\ast$ defined on the relative 
cohomology group is not a monomorphism.
\end{proof}

\section{Hyperelliptic Components}
\label{s:hyperelliptic:components}

\subsection{The Teichm\"uller Geodesic Flow}

The group $GL(2,\mathbb{R})_{+}$ acts naturally on the moduli space 
of quadratic differentials; this action preserves the natural stratification 
of the moduli space.

The action of the diagonal subgroup of $SL(2,\mathbb{R})$ on the moduli space 
of quadratic differentials can be naturally identified with the Teichm\"uller 
geodesic flow on the moduli space of curves for the Teichm\"uller metric. 

Many important results in the theory of interval exchange maps, of measured 
foliations and of dynamics on translation surfaces are based on the following 
fundamental Theorem, proved independently by 
H.Masur~\cite{Masur:82} and by W.Veech~\cite{Veech:82}:

\begin{NoNumberTheorem}[H.~Masur; W.~Veech]
The Teichm\"uller geodesic flow acts ergodically on every 
connected component of every (normalized) stratum of the moduli space 
of quadratic differentials; the corresponding 
invariant measure is a finite Lebesque equivalent measure.
\end{NoNumberTheorem}

\begin{Proposition}
\label{prop:releant}
The action of $GL(2,\mathbb{R})_{+}$ on the strata 
commutes with the mapping of Construction~\ref{constr:general}.
\end{Proposition}

\begin{proof}
In the proof of Lemma~\ref{lm:mapping:is:not:degenerate}, we use some 
coordinates to linearize the mapping. It is easy to see 
that in these charts, the linear action commute with the mapping.
\end{proof}

\subsection{Hyperelliptic Components}

Let us apply the Construction~\ref{constr:general} in the
following particular case.  Consider a meromorphic quadratic
differential $\psi$ on $\CP$ having the singularity pattern
$(2(g-k)-3,2k+1,-1^{2g+2})$, where $k \geq -1$, $g \geq 1$ and
$g-k \geq 2$. Consider a ramified double covering $\pi$ over $\CP$
having ramification points over $2g+2$ poles of $\psi$, and no
other ramification points. We obtain a hyperelliptic Riemann
surface $\tilde{M}$ of genus $g$ with a quadratic differential
$\pi^{*}\psi$ on it. By Lemma~\ref{lm:change:of:degree} the
induced quadratic differential $\pi^{*}\psi$ has the singularity
pattern $(2(g-k)-3,2(g-k)-3,2k+1,2k+1)$.
Construction~\ref{constr:general} gives us a local mapping
$$
\Q(2(g-k)-3,2k+1,-1^{2g+2}) \to \Q(2(g-k)-3,2(g-k)-3,2k+1,2k+1),
$$
where $k \geq -1, g \geq 1$ and $g-k \geq 2$. Computing the
dimensions of the strata we get
\begin{align*}
&\dim_{\C}\Q(2(g-k)-3,2k+1,-1^{2g+2})=2\cdot 0 + (2g + 4) - 2=2g+2\\
&\dim_{\C}\Q(2(g-k)-3,2(g-k)-3,2k+1,2k+1)=2 g + 4 -2 =2g+2
\end{align*}
Thus the dimensions of the strata coincide. By
Lemma~\ref{lm:mapping:is:not:degenerate} the mapping is
non-degenerate.

When $2k+1\neq -1$ the mapping
$$
\Q(2(g-k)-3,2k+1,-1^{2g+2})
\to \Q(2(g-k)-3,2(g-k)-3,2k+1,2k+1)
$$
is globally defined.
Since the stratum
$\Q(2(g-k)-3,2k+1,-1^{2g+2})$ is connected, as any
other stratum on $\CP$ by Proposition in
section~\ref{s:covering:construction}, see above, the
image of the mapping in the stratum $\Q(2(g-k)-3,2(g-k)-3,2k+1,2k+1)$ is also
connected. When $2k+1= -1$ the mapping has
$2g+3$ branches corresponding to the choice of the simple pole of
$\psi$ where we do not have a ramification of the covering.
However, since we can deform the positions of the zero and the
poles of $\psi$ on $\CP$ arbitrarily (avoiding collapses, of
course), the intersection of the image is non empty,
so the union of the images of the
mapping $ \Q(2g-1,-1^{2g+3}) \to \Q(2g-1,2g-1,-1,-1) $ is
connected. 

Since the dimension of the strata coincide, and the 
mapping is non-degenerate, we obtain an open set on the stratum 
$\Q(2(g-k)-3,2(g-k)-3,2k+1,2k+1)$. By Proposition~\ref{prop:releant}, 
the action of the geodesic flow is relevant, thus by ergodicity of 
this flow, the image of the mapping
$$
\Q(2(g-k)-3,2k+1,-1^{2g+2}) \to
\Q(2(g-k)-3,2(g-k)-3,2k+1,2k+1)
$$
gives us a full measure set in the corresponding connected component of 
the stratum
$$
\Q(2(g-k)-3,2(g-k)-3,2k+1,2k+1).
$$
Thus we obtain a connected component of these stratum.

Similarly to the previous case we can easily check coincidence of
the dimensions of the strata
 $$
\Q(2(g-k)-3,2k,-1^{2g+1}) \to \Q(2(g-k)-3,2(g-k)-3,4k+2),
 $$
with $k \geq 0$, $g \geq 1$ and $g-k \geq 1$ and
 $$
 \Q(2g-2k-4,2k,-1^{2g}) \to \Q(4(g-k)-6,4k+2)
 $$
with $k \geq 0$, $g \geq 2$ and $g-k \geq 2$.

The images of these mappings give us connected components in the
strata
 $$
 \Q(2(g-k)-3,2(g-k)-3,4k+2).
 $$
and
 $$
 \Q(4(g-k)-6,4k+2).
 $$

\begin{Definition}
\label{def:main:theo}
The connected components constructed above are called the {\it
hyperelliptic components} and are denoted by:
\begin{enumerate}
\item $\Q(2(g-k)-3,2k+1,-1^{2g+2}) \to \Q^{hyp}(2(g-k)-3,2(g-k)-3,2k+1,2k+1),$\newline
where $k \geq -1, \ g \geq 1, \ g-k \geq 2$.
The corresponding double covering has ramification points over
$2g+2$ poles of meromorphic quadratic differential on $\CP$.
\item $\Q(2(g-k)-3,2k,-1^{2g+1}) \to \Q^{hyp}(2(g-k)-3,2(g-k)-3,4k+2)$,\newline
where $k \geq 0$, $g \geq 1$ and $g-k \geq 1$.
The corresponding double covering has ramification points over
$2g+1$ poles and over the zero of degree $2k$
of meromorphic quadratic differential on $\CP$.
\item $\Q(2g-2k-4,2k,-1^{2g}) \to \Q^{hyp}(4(g-k)-6,4k+2)$,\newline
where $k \geq 0$, $g \geq 2$ and $g-k \geq 2$.
The corresponding double covering has ramification points over
all singularities of the quadratic differential on $\CP$.
\end{enumerate}
\end{Definition}

\begin{Remark}
The connected component $\Q^{hyp}(2(g-k)-3,2(g-k)-3,2k+1,2k+1)$ was
first noticed by M.~Kontsevich.
\end{Remark}

\begin{Remark}
For the mapping $$\Q(2g-2k-4,2k,-1^{2g}) \to \Q(4(g-k)-6,4k+2)$$one can see
that the resulting quadratic differential is not the global square
of an Abelian differential in the following way. Take a path $\gamma$ on $\CP$
around the zero of order $2k$ and a pole of the corresponding quadratic
differential. The monodromy of
the covering along this path is trivial so it can be lifted on a path
$\tilde{\gamma}$ in $M$. The holonomy of the flat structure
along the path $\tilde{\gamma}$ is non trivial,
we can compute how the tangent vector to the corresponding path
$\tilde{\gamma}$ on $M$ turns in the flat structure defined by
$\psi$ (see Figure~\ref{fig:index}). The counterclockwise
direction is chosen as a direction of the positive turn.

\begin{figure}[htbp]
\begin{center}
\psfrag{p1}{\scriptsize $\tilde{P_1}$}  \psfrag{p2}{\scriptsize $\tilde{P_2}$}
\psfrag{a1}{\scriptsize $A_1$}  \psfrag{a5}{\scriptsize $A_5$}  
\psfrag{a2}{\scriptsize $A_2$}  \psfrag{a6}{\scriptsize $A_6$}  
\psfrag{a3}{\scriptsize $A_3$}  \psfrag{a7}{\scriptsize $A_7$}  
\psfrag{a4}{\scriptsize $A_4$}  \psfrag{a8}{\scriptsize $A_8$}  

\includegraphics[width=8cm]{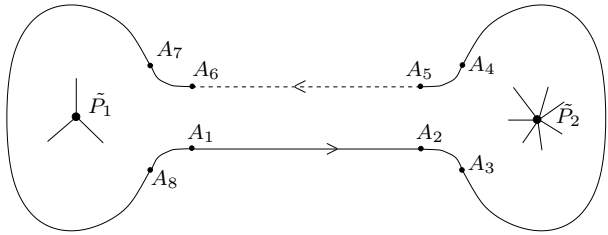}
\end{center}
\caption{
\label{fig:index}
A path $\tilde{\gamma}$ on the covering surface obtained by lifting of the path
$\gamma$. Since $\tilde{P}_1, \tilde{P}_2\in \tilde{M}$ are ramification
points, the segments $A_1 A_2$ and $A_6 A_5$ are located at
the ``different sheets'' of the covering $\tilde{M}\to M$, and
not nearby.
}
\end{figure}

We get some angle $\phi$ following the part $A_1 A_2$ of the path
which goes from one singularity to another. Then we make a turn by
$-\pi/2$ going along $A_2 A_3$. Turning around the singularity
(the path $A_3 A_4$) we get the angle $(k_2+2)\pi$, which is
followed by another turn by $-\pi/2$ now along $A_4 A_5$. The path
$A_5 A_6$ gives the turn by $-\phi$, which is followed by another
$-\pi/2$ along $A_6 A_7$. Turning around singularity along $A_7
A_8$ we get the angle $(k_1+2)\pi$, and the loop is competed by the
path $A_8 A_1$ giving one more turn by $-\pi/2$. All together this
gives $(k_1+k_2+2) \cdot \pi$.

In our special case, the resulting angle for the corresponding 
quadratic differential is $(2k+1)\cdot\pi$, thus the corresponding
measured foliation is not oriented and consequently, the quadratic 
differential $\psi$ on $M$ is not the 
global square of an Abelian differential. The image of application 
above belong in strata $\Q(4(g-k)-6,4k+2)$ and the application above 
is well defined.
\end{Remark}

\begin{Remark}
\label{particular:map}
In the each of the following cases
\begin{align}
&\Q(-1,-1,-1,-1,0,0) \to \Q(2,2)
\notag \\
\label{eq:particular:map}
&\Q(-1,-1,-1,-1,0) \to \Q(-1,-1,2)
\\
&\Q(-1,-1,-1,-1) \to \Q(-1,-1,-1,-1)
\notag
\end{align}
the image stratum also has the same dimension as the original one.
All the corresponding ramified coverings have even degree $2d\ge 2$.

In the first case the covering has $d$ ramification points of degree
$2$ over each of four simple poles of the meromorphic quadratic differential
on $\CP$; a single ramification point of degree
$2$ over each of two marked points, and no other ramification points.

In the second case the covering has $d-1$ ramification points of degree
$2$ over one of the poles; $d$ ramification points of degree
$2$ over each of the remaining three simple poles of the meromorphic quadratic differential
on $\CP$; a single ramification point of degree
$2$ over the marked point, and no other ramification points.

In the third case the covering has $d-1$ ramification points of degree
$2$ over each of two poles; $d$ ramification points of degree
$2$ over each of remaining two simple poles of the meromorphic
quadratic differential on $\CP$, and no other ramification points.

The strata $\Q(-1,-1,2)$ and $\Q(2,2)$ are in the list of hyperelliptic
components in the definition above; the corresponding surfaces
of genus $1$ and $2$ are respectively elliptic and hyperelliptic. We
show in section~\ref{s:qd:versa:foliations}, Lemma~\ref{lm:p:1} that
theses strata are connected, which implies that the corresponding
mappings have no interest for our purposes.
The stratum  $\Q(-1,-1,-1,-1)$ corresponds to quadratic differentials
on $\CP$ so it is connected by Proposition of Kontsevich.
Thus the third mapping is of no interest for us neither.
\end{Remark}

The Theorem below asserts that there are no other connected
components which can be obtained using a similar construction.

\begin{Theorem}
\label{thm:main}
Let $\Q(k_1,\dots,k_n)$ be a stratum in the moduli space of
meromorphic quadratic differentials and let $\pi : \tilde{M}
\rightarrow M$ be a covering of finite degree $d > 1$. Consider
the mapping
 $$
  \Q(k_1, \dots, k_n) \to \Q(\tilde{k}_1,\dots,\tilde{k}_m)
 $$
induced by the covering $\pi$ (see
Construction~\ref{constr:general}). Suppose that the image stratum
is not a stratum of squares of Abelian differentials, and suppose
that the mapping is neither
of one of the three types corresponding to
hyperelliptic components nor of one of the three exceptional
types~\eqref{eq:particular:map}. Then
$$
\dim_{\C} \Q(k_1,\dots,k_n) < \dim_{\C}
\Q(\tilde{k}_1,\dots,\tilde{k}_m)
$$
\end{Theorem}
\begin{proof}

We introduce the integer parameters $d,n,m,p,r$ responsible
for the topological type of the pair: (covering $\pi$, quadratic
differential $\psi$).
Let $d$ denote the degree of the covering, 
$n$ the number of true zeros which are critical values,
$m$ the number of marked points, 
$p$ the number of simple poles,
$r$ the number of singularities of $\psi$ which are regular 
values for the covering $\pi$.
See below for an explicit example and details.

In the first part of
the proof we derive from the relation $\dim\Q=\dim\tilde{\Q}$
elementary inequality~\eqref{weakened:inequality} for the positive
integers $d,n,m,p,r$. The inequality has not so many solutions.
In Lemma~\ref{lm:d:ge:3} we show that solutions with $d\ge 3$,
where $d$ stands for the number of branches of the covering $\pi$, do
not correspond to any nontrivial mappings $\Q\to\tilde{\Q}$.
In Lemma~\ref{lm:d:2} we show that the only two-fold coverings
$\pi$ which give rise to the solution of the equation
$\dim\Q=\dim\tilde\Q$ are exactly those which correspond to
hyperelliptic components.

\begin{figure}[htbp]
\begin{center}
\psfrag{pi}{$\pi$}
\psfrag{m}{$M$}
\psfrag{mt}{$\tilde M$}

\psfrag{t}{$\underbrace{{\scriptstyle P_1,\dots,P_n}}_{\textrm{True zeros}}$}
\psfrag{ma}{$\underbrace{{\scriptstyle P_{n+1},\dots,P_{n+m}}}_{\textrm{Marked points}}$}
\psfrag{s}{$\underbrace{{\scriptstyle P_{n+m+1},\dots,P_{n+m+p}}}_{\textrm{Simple Poles}}$}
\psfrag{n}{$\underbrace{{\scriptstyle P_{n+m+p+1},\dots,P_{n+m+p+r}}}_{\textrm{Non critical
values of the covering $\pi$}}$}

\includegraphics[width=12cm]{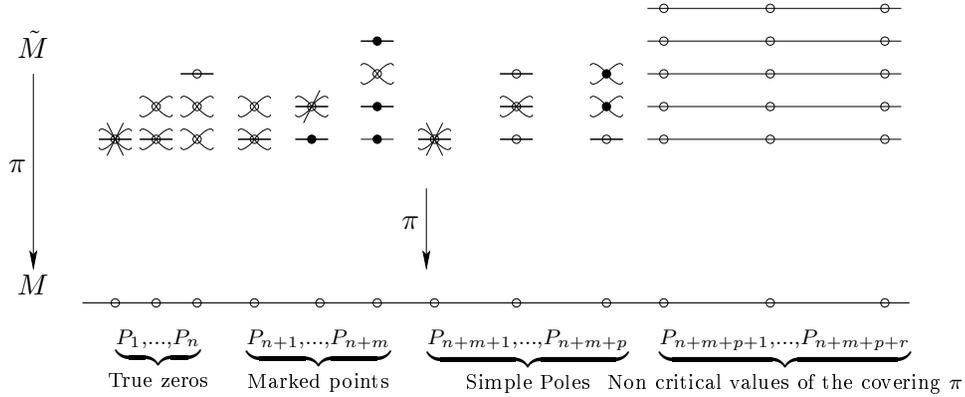}
\end{center}

\vskip 5mm
\caption{
\label{fig:covering}
This example presents a covering of degree $5$ with $21$
ramification points. There are $9$ critical values of $\pi$.
Here we have $n=m=p=r=3$.
The points $\tilde{P} \in \pi^{-1}(P_i)$ which correspond to regular
points of $\pi^\ast\psi$ are painted in black. By convention we
{\it do not} mark them.}
\end{figure}

Let us first introduce some notations. By $(M,\psi)$ we denote the
initial Riemann surface and a quadratic differential on it. We
consider a ramified covering $\pi: \tilde{M}\to M$ of degree $d$.
We denote by $P_1, \dots, P_n \in M$ be the critical values corresponding
to the {\it true} zeros of $\psi$, i.e. to the zeros of degrees
$k_i>0$. Let $P_{n+1}, \dots, P_{n+m} \in M$ be the critical
values corresponding to the marked points (or, in the other words,
to the {\it fake} zeros of $\psi$).  Let $P_{n+m+1}, \dots,
P_{n+m+p} \in M$ be the critical values corresponding to the
simple poles of $\psi$. Let $P_{n+m+p+1}, \dots, P_{n+m+p+r} \in M$
be the remaining true zeros and simple poles of $\psi$ which
correspond to non-critical values of $\pi$. (See an example at
Figure~\ref{fig:covering}).


We consider only those ramified coverings $\pi$ for which the
induced quadratic differential $\tilde{\psi}=\pi^\ast\psi$ on
$\tilde{M}$ is not globally a square of a holomorphic 1-form. We
do not exclude from consideration the situation when the covering
$\pi$ is actually a regular covering.


Let us discuss now the location of singularities of the quadratic
differential $\tilde{\psi}=\pi^\ast\psi$ on $\tilde{M}$. By
convention we do not endow $\tilde{\psi}$ with any {\it fake}
zeros: all singularities of $\tilde{\psi}$ are either {\it true}
zeros, or simple poles.


The computation of Lemma~\ref{lm:change:of:degree} shows that
$\tilde{\psi}$ may have a singularity at a point
$\tilde{P}\in\tilde{M}$ either when $\psi$ has a singularity at
$\pi(\tilde{P})$ or when $\tilde{P}$ is a ramification point of
$\pi$. Thus all zeros and poles of $\tilde{\psi}$ are located at
the preimage of the points $P_1, \dots, P_{n+m+p+r}$. Let us
specify them precisely. (See also Figure~\ref{fig:covering} where
singularities are denoted by white bullets.)


At any preimage of any of the zeros $P_1,\dots, P_n$ of $\psi$ we
have a true zero of $\tilde{\psi}$, see
Lemma~\ref{lm:change:of:degree}. We have a zero of $\tilde{\psi}$
at a preimage $\tilde{P}\in\pi^{-1}(P_i)$ of a marked point $P_i$,
$n+1\le i \le n+m$, if and only if $\tilde{P}$ is a ramification
point of the covering $\pi$, $e_{\pi}(\tilde{P})\neq 1$. We have a
pole or a true zero of $\tilde{\psi}$ at a preimage
$\tilde{P}\in\pi^{-1}(P_i)$ of a simple pole $P_i$, $n+m+1\le i
\le n+m+p$, if and only if the ramification index
$e_{\pi}(\tilde{P})$ of the covering $\pi$ at $\tilde{P}$ is
different from two, $e_{\pi}(\tilde{P})\neq 2$. At any preimage of
any singularity $P_{n+m+p+1},\dots, P_{n+m+p+r}$ of $\psi$ we
have a simple pole or a true zero of $\tilde{\psi}$. (See also an
example presented at Figure~\ref{fig:covering}.)


Now everything is ready for the computation of dimensions of the
strata of meromorphic quadratic differentials $\Q$ and
$\tilde{\Q}$   corresponding to $\psi$ and $\tilde{\psi}$.

Let $g$ be the genus of the Riemann surface $M$. We have
\begin{equation}
\label{eq:dim:X}
 \dim\Q = 2g + n+m+p+r -2
\end{equation}
   
Denoting by $\tilde{g}$ be the genus of the Riemann surface
$\tilde{M}$ we get
   
\begin{align}
 \dim\tilde{\Q} &= 2\tilde{g} +
 \text{number of singularities of } \tilde{\psi} -2 =
\notag
\\
 &=  2\tilde{g} +
 \sum_{i=1}^{n}\text{number of preimage $\tilde{P}$ of } \pi^{-1}(P_i)
 \quad +\notag\\& \qquad +
 \sum_{i=n+1}^{n+m}\text{number of preimage $\tilde{P}$ of
 $\pi^{-1}(P_i)$ of index } e_{\pi}(\tilde{P})\neq 1
 \quad +\notag\\& \qquad +
 \sum_{i=n+m+1}^{n+m+p}\text{number of preimage $\tilde{P}$ of
 $\pi^{-1}(P_i)$ of index }e_{\pi}(\tilde{P})\neq 2
 \quad +\notag\\& \qquad +
 \sum_{i=n+m+p+1}^{n+m+p+r}\text{number of preimage $\tilde{P}$ of }
 \pi^{-1}(P_i)\quad - \quad 2 \quad =\notag
\\
\label{eq:dim:tildeQ}
 &=  (2\tilde{g}-2) +d\cdot r +
 \sum_{i=1}^{n+m+p}\text{number of preimage $\tilde{P}$ in }
 \pi^{-1}(P_i)\quad -
 \\&\qquad-
 \sum_{i=n+1}^{n+m}\text{number of preimage $\tilde{P}$ in
 $\pi^{-1}(P_i)$ of ramification index one}\ -
 \notag\\&\qquad-
 \sum_{i=n+m+1}^{n+m+p}\text{number of preimage $\tilde{P}$ in
 $\pi^{-1}(P_i)$  of ramification index two}\notag
\end{align}
   
Computing the genus $\tilde{g}$ of the Riemann surface
$\tilde{M}$ by Riemann--Hurwitz formula
   
\begin{align*}
2\tilde{g}-2 & = d\cdot (2g-2)+
 \sum_{\text{ramification points $\tilde{P}$ of $\pi$}} (e_{\pi}(P)-1)
 =
\\
 & = d\cdot (2g-2)+
 \sum_{\text{critical values $P_i$ of $\pi$}} (d - \text{number of preimage})
\\
 & = d\cdot (2g-2)+ d\cdot(n+m+p) -
 \sum_{i=1}^{n+m+p}\text{number of preimage $\tilde{P}$ in }
 \pi^{-1}(P_i)
\end{align*}
and substituting $2\tilde{g}-2$ in the
formula~\eqref{eq:dim:tildeQ} by the latter expression we get the
following answer:
\begin{align}
\label{dim:tildeQ}
 &\dim\tilde{\Q} =  d\cdot (2g+n+m+p+r-2) -
 \\
 &\qquad-
 \sum_{i=n+1}^{n+m}\text{number of preimage $\tilde{P}$ in
 $\pi^{-1}(P_i)$ of ramification index one}\ -
 \notag\\&\qquad-
 \sum_{i=n+m+1}^{n+m+p}\text{number of preimage } \tilde{P}\in\pi^{-1}(P_i)
 \text{ of ramification index two}\notag
\end{align}
Comparing the dimensions $\dim\Q$ and $\dim\tilde{Q}$ given by
equations~\eqref{eq:dim:X} and~\eqref{dim:tildeQ} we see that the
equation
$$
  \dim\Q = \dim\tilde{Q}
$$
is equivalent to the following one
\begin{multline}
\label{eq:dimQ:dim:tildeQ}
 (d-1)\cdot (2g+n+m+p+r-2) =
 \\
 =
 \sum_{i=n+1}^{n+m}\text{number of preimage $\tilde{P}$ in
 $\pi^{-1}(P_i)$ of ramification index one}\ +
 \\ +
 \sum_{i=n+m+1}^{n+m+p}\text{number of preimage } \tilde{P}\in\pi^{-1}(P_i)
 \text{ of ramification index two}
\end{multline}

By convention $P_{n+1}, \dots, P_{n+m}\in M$ are critical values
of the $d$-fold covering $\pi$. Hence, among preimage of any
$\tilde{P}$ in $\pi^{-1}(P_i)$, $n+1\le i \le n+m$ there is at
least one preimage with ramification index different from one.
This implies that the number of preimage $\tilde{P}$ in
$\pi^{-1}(P_i)$ of ramification index one is at most $d-2$.

The covering $\pi:\tilde{M}\to M$ is of order $d$. Hence, for any
point $P_i\in M$ the number of preimage
$\tilde{P}\in\pi^{-1}(P_i)$ of ramification index two has the
following obvious upper bound
\begin{multline}
\label{eq:preimage:of:pole}
 \text{number of preimage } \tilde{P} \text{ in } \pi^{-1}(P_i)
 \text{ of ramification index two}  \\ \le
 \begin{cases}
    d/2 & \text{ when $d$ is even }\\
    (d-1)/2 & \text{ when $d$ is odd }
 \end{cases}
 \end{multline}
Using these two obvious bounds we derive the following inequality
from equation~\eqref{eq:dimQ:dim:tildeQ}
\begin{multline}
\label{basic:inequality}
 (d-1)\cdot (2g+n+m+p+r-2)\le
\begin{cases}
     m\cdot(d-2) + p\cdot d/2 & \text{ when $d$ is even }\\
     m\cdot(d-2) + p\cdot(d-1)/2 & \text{ when $d$ is odd }
 \end{cases}
\end{multline}
and its weakened version
\begin{equation}
\label{weakened:inequality}
 (d-1)\cdot (2g-2+n+r) \le p\cdot(1-d/2) - m
\end{equation}

In what follows we consider the nonnegative integer solutions of
inequalities~\eqref{basic:inequality}
and~\eqref{weakened:inequality}. There are not so many of them
since for most types of the covering $\pi$ the dimension
$\dim\tilde{Q}$ is greater than $\dim\Q$ and the parameters
$d,n,m,p,r$ do not obey neither
equation~\eqref{eq:dimQ:dim:tildeQ} nor even
inequalities~\eqref{basic:inequality}
and~\eqref{weakened:inequality}. In Lemma~\ref{lm:d:ge:3} we show
that the solutions of inequalities with $d\ge 3$ do not correspond
to any actual ramified coverings $\pi$. In Lemma~\ref{lm:d:2} we
study solutions of inequalities corresponding to $d=2$, and show
that the two-fold coverings $\pi$ giving solutions of
equation~\eqref{eq:dimQ:dim:tildeQ} are exactly those which
corresponding to hyperelliptic components.

\begin{Lemma}
\label{lm:d:ge:3}
If $d\ge 3$ then $ \dim\Q < \dim\tilde{\Q}$.
\end{Lemma}
\begin{proof}[Proof of Lemma~\ref{lm:d:ge:3}]

Suppose that dimensions of $\Q$ and $\tilde{\Q}$ coincide. We can apply
the inequalities above.

First of all, we are going to show that genus $g$ of the underlying surface
$M$ must be equal to $0$.
Suppose that $g \geq 1$, so inequality~\eqref{weakened:inequality} gives
$(d-1)(n+r) \leq -p/2 -m$. If $p > 0$ or $m > 0$ then $n+r$ must be negative,
which is impossible. If $p=0$ and $m=0$
then $n=r=p=m=0$ and we obtain the stratum $\Q(\emptyset)$ which is empty.
Thus the assumption $g \geq 1$ leads to a contradiction and we must have
$g=0$ which means that $M \simeq \CP$.

Now suppose that $m > 0$. Since the right-hand side of
inequality~\eqref{weakened:inequality} is strictly negative,
and $d-1>0$,
the expression $(2g-2+n+r)$ should be strictly negative. Since
$g=0$ this implies that
$$
n+r \leq 1
$$

Now suppose that $m=0$. Using inequality~\eqref{weakened:inequality}, we get
$$(d-1)\cdot(n+r-2) \leq p\cdot(1-d/2) \le -p/2$$
If $p=0$ then $n+r \leq 2$ so $r+p \leq 2$ which is impossible
since a meromorphic quadratic differential on $\CP$ has
at least four poles. Thus $p > 0$ and we again get
$n+r \leq 1$.

The above remarks show that we can restrict our considerations
to the following three (overlapping) cases:
\begin{displaymath}
\begin{array}{ccc}
\begin{cases}
n+r = 0 \\
m = 0
\end{cases}  &
\begin{cases}
n+r=1 \\
m \geq 0
\end{cases}  &
\begin{cases}
n+r=0 \\
m > 0
\end{cases}
\end{array}
\end{displaymath}
To finish the proof of Lemma~\ref{lm:d:ge:3} we are going to show that
all possible mappings with $d\ge 3$ satisfying~\eqref{eq:dimQ:dim:tildeQ}
are listed in~\eqref{eq:particular:map}. We shall use the following
obvious remark which is valid for any ramified covering
$\pi: \tilde{M}\to M$ of degree $d$ and any point $P\in M$:
\begin{equation}
\label{eq:fiber}
\sum_{\tilde{P}\in\pi^{-1}(P)} e_\pi(\tilde{P}) = d
\end{equation}

\begin{Caseone}
Since $n=r=m=0$ and genus $g=0$, we get a meromorphic quadratic
differential with four simple poles $P_1, \dots, P_4$
on $\CP$, $p=4$. Let $t_i$ be the
number of points in the fiber
$\pi^{-1}(P_i)$ ($i=1,\dots,4$) of ramification index two.

Taking the sum of expressions~\eqref{eq:fiber} over $P_i$, $i=1,\dots,4$, we
get
\begin{multline*}
4d=\sum_{i=1}^4\sum_{\tilde{P}\in\pi^{-1}(P_i)} e_\pi(\tilde{P})
= \sum_{i=1}^4\Big(
   \sum_{\substack{\tilde{P}\in\pi^{-1}(P_i)\\
         e_\pi(\tilde{P}) \not = 2 }}
      e_\pi(\tilde{P}) +
   \sum_{\substack{\tilde{P}\in\pi^{-1}(P_i)\\
         e_\pi(\tilde{P}) = 2 }}
      e_\pi(\tilde{P})
\Big)=
\\
= \sum_{i=1}^4
   \sum_{\substack{\tilde{P}\in\pi^{-1}(P_i)\\
         e_\pi(\tilde{P}) \not = 2 }}
      e_\pi(\tilde{P})
 + 2\sum_{i=1}^4 t_i
\end{multline*}
On the other hand equation~\eqref{eq:dimQ:dim:tildeQ} gives
\begin{equation*}
(d-1)\cdot(4-2) = \sum_{i=1}^4 t_i
\end{equation*}
Note that, as $m=0$, the term corresponding to  
ramification index $1$ in equation~\eqref{eq:dimQ:dim:tildeQ}
does not appear. \\
Thus we obtain
$$
\sum_{i=1}^4
   \sum_{\substack{\tilde{P}\in\pi^{-1}(P_i)\\
         e_\pi(\tilde{P}) \not = 2 }}
      e_\pi(\tilde{P}) = 4
$$
This means that either there are two ramification points
$\tilde{P_1},\tilde{P_2}$ of ramification index different from
$2$, and then $e_\pi(\tilde{P_1})=1, e_\pi(\tilde{P_2})=3$, or
there are four ramification points $\tilde{P_i},i=1,\dots,4$ of
ramification index different from $2$, and then
$e_\pi(\tilde{P_i})=1$.

The first solution suggests the map of
the moduli spaces
$$
\Q(-1,-1,-1,-1) \longrightarrow \Q(-1,1).
$$
But the stratum $\Q(-1,1)$ is empty (see above). The second
solution corresponds  to the map $\Q(-1,-1,-1,-1) \longrightarrow
\Q(-1,-1,-1,-1)$ which is one of the
maps~\eqref{eq:particular:map}, and hence has no interest for us.
\end{Caseone}

\begin{Casetwo}
In this case either $n=1$ or $r=1$; denote by $k$ the order of
corresponding singularity of $\psi$. The singularity is a true
zero if $k > 0$ or a pole if $k=-1$; by convention on notations
$n,r$ it cannot be a marked point. With this notation $p=4+k$. We
derive from~\eqref{weakened:inequality} the following inequality:
$-(d-1) \leq (4+k)\cdot(1-d/2)-m$ which gives $d \leq 3-k/2-m$. If
$m \geq 1$, as $d \geq 3$, we obtain $k \leq -2$ which is
forbidden since we consider only simple poles. Thus we must have
$m=0$ and so $k=-1$. This implies that $n=0$ and $r=1$.

By inserting in inequality~\eqref{weakened:inequality} we obtain
$d=3 \textrm{ or } 4$. The stronger
inequality~\eqref{basic:inequality} eliminates the solution $d=3$,
so finally we get the solution $d=4$, where the values of the
other parameters are as follows: $g=n=m=0$; $r=1$, and the
corresponding point is a simple pole; there are three more simple
poles, so $r=3$.

However, this solution gives the map $\Q(-1,-1,-1,-1) \rightarrow
\Q(-1,-1,-1,-1)$ of the moduli spaces from the exceptional
list~\eqref{eq:particular:map}, and hence has no interest for us.
\end{Casetwo}

\begin{Casethree}
Since $n=r=0$ and genus $g=0$, we get a meromorphic quadratic
differential with $p=4$ simple poles on $\CP$ and with $m>0$
marked points. Using inequality~\eqref{basic:inequality} we derive
that if $d$ is odd, we have $m=0$ which is a contradiction. So $d$ 
must be even and again by using this inequality with $d$ even, we 
deduced that $1 \leq m \leq 2$.

If $m=1$, then we have a single marked point $P_0\in M$. We derive
the following relation from equation~\eqref{eq:dimQ:dim:tildeQ}:
\begin{multline*}
 (d-1)\cdot(1+4-2) =
 \\
 =
 \sum_{\textrm{poles $P_i$}}\text{number of preimage $\tilde{P}$ in
 $\pi^{-1}(P_i)$ of ramification index two}\ +
 \\ +
 \text{number of preimage } \tilde{P}\in\pi^{-1}(P_0)
 \text{ of ramification index one}
\end{multline*}
Let $t$ be the number of preimage $\tilde{P}\in\pi^{-1}(P_0)$ of
ramification index one over the marked point $P_0$. By
Convention, at least
one preimage of a marked point has ramification index greater
then one, thus~\eqref{eq:fiber} implies $t \leq d-2$.
Taking into consideration~\eqref{eq:preimage:of:pole}
we obtain from the equality above the following inequality
$$
3\cdot(d-1) \leq 2\cdot d + t
$$
and hence we have
$$
d-3 \leq t \leq d-2
$$
It is not difficult to see that if $t=d-3$ we obtain the  following map
$$
\Q(-1,-1,-1,-1,0) \to \Q(4)
$$
but $\Q(4)$ is empty.

If $t=d-2$ the equality above implies that
$$
\sum_{\textrm{$P_i$ poles}}\text{number of preimages $\tilde{P}$ in
 $\pi^{-1}(P_i)$ of ramification index two}\ = 2\cdot d -1
$$
and taking in consideration~\eqref{eq:fiber} we get
$$
\sum_{i=1}^4
   \sum_{\substack{\tilde{P}\in\pi^{-1}(P_i)\\
         e_\pi(\tilde{P}) \not = 2 }}
      e_\pi(\tilde{P}) = 2
$$
where $P_1, \dots, P_4$ are the simple poles.
It is not difficult to see that in this case we obtain the map
$$
\Q(-1,-1,-1,-1,0) \to \Q(-1,-1,2)
$$
which belongs to the list~\eqref{eq:particular:map}
and hence does not interest us.

By similar arguments we conclude that
if $m=2$ we get the map
$$
\Q(-1,-1,-1,-1,0,0) \to \Q(2,2)
$$
which belongs to the list~\eqref{eq:particular:map}
and hence does not interest us.
\end{Casethree}
     %
Lemma~\ref{lm:d:ge:3} is proved.
\end{proof}

It remains to consider the two-fold coverings in order to complete
the proof of Theorem~\ref{thm:main}.

\begin{Lemma}
\label{lm:d:2}
If $d=2$ then the strata obey the relation
$$
\dim\Q = \dim\tilde{\Q}
$$
only in the cases listed in Theorem~\ref{thm:main} and in
exceptional list~\eqref{eq:particular:map}.
\end{Lemma}
\begin{proof}[Proof of Lemma~\ref{lm:d:2}]
We suppose that the dimensions coincide. As the degree of the covering
is two, there are no ramification points in $\tilde{M}$ of index one over
marked points $P_{n+1},\dots,P_{n+m}$ on $M$, and
there is exactly one point
over each pole
$P_{n+m+1},\dots,P_{n+m+r}$ in the set of critical
values of $\pi$. Thus the
equation~\eqref{eq:dimQ:dim:tildeQ} transforms into the following one
$$2g+n+m+p+r-2 = 0 + p = p$$
$$2g+n+m+r-2=0$$
This equation shows that the genus $g$ of the underlying surface
might be either $1$ or $0$. If $g=1$ we get $n=m=r=0$ and hence
the meromorphic quadratic differential on a surface of genus
$g=1$ does not have either zeros nor marked points.
Hence it does not have simple poles neither, which implies that
it is the square of an Abelian differential.

Thus $g=0$ and $n+m+r=2$. We consider separately the following three
cases:
\begin{displaymath}
\begin{array}{ccc}
\begin{cases}
n+m=0 \\
r=2
\end{cases}  &
\begin{cases}
n+m=1 \\
r=1
\end{cases}  &
\begin{cases}
n+m=2 \\
r=0
\end{cases}
\end{array}
\end{displaymath}
The first two cases give us either the first two maps from
Definition~\ref{def:main:theo} or the exceptional 
map $\Q(-1,-1,-1,-1)\to \Q(-1,-1,-1,-1)$ from the list~\eqref{eq:particular:map}.

The third case gives the mapping (when  it is defined)
\begin{equation}
\label{eq:third:case}
 \Q(k_1,k_2,-1^{k_1+k_2+4}) \rightarrow \Q(2k_1+2,2k_2+2)
\end{equation}
with $k_i \geq 0$.

A ramified double covering over $\CP$ has even number of
ramification points. Since in our case we have ramification
points over the two singularities of degrees $k_1$ and $k_2$ and
over all simple poles, this implies that the number $k_1+k_2+6$ is
even, and hence $k_1$ and $k_2$ has the same parity. When both of
them are odd, one can recognize in~\eqref{eq:third:case} the
canonical ramified double covering described in
Construction~\ref{constr:canonical:2:covering}. Thus in this case
the resulting quadratic differential is the global square of an
Abelian differential, and this case does not interest us.

When both $k_1, k_2 \ge 0$ are even we obtain the map of the third
type from Definition~\ref{def:main:theo}.
This completes the proof of Lemma~\ref{lm:d:2}.
\end{proof}

Theorem~\ref{thm:main} now follows immediately from
Lemmas~\ref{lm:d:ge:3} and~\ref{lm:d:2}.

\end{proof}

\section{Quadratic Differentials Versa Flat Structures}
\label{s:qd:versa:foliations}

In this section we present a well-know relation between quadratic
differentials and flat structures on Riemann surfaces.
We use this relation to prove that the
strata, that possess a hyperelliptic connected component,
are not connected except several particular cases in low genera.

A flat surface with cone type singularities is a surface which
possesses locally the geometry of a standard cone. We can define
it by a flat Riemannian metric with specific isolated
singularities. The standard cone possesses a unique invariant: it
is the angle at the vertex. Here we consider only {\it
half-translation} flat surfaces: parallel transport of a tangent
vector along any closed path either brings the vector $\vec{v}$
back to itself or brings it to the centrally-symmetric vector
$-\vec{v}$. This implies that the cone angle at any singularity of
the metric is an integer multiple of $\pi$.

Let $\psi$ be a meromorphic quadratic differential on a Riemann
surface $M^2_g$. Then it is possible to construct an atlas on $M
\backslash \{$singularities$\}$ such that $\psi=dz^2$ in any
coordinate chart. As $dz^2=dw^2$ implies $z=\pm w +const$, we see
that the charts of the atlas are identified either by a
translation or by a translation composed with a central symmetry.
Thus a meromorphic quadratic differential $\psi$ induces a
half-translation flat structure on $M \backslash
\{$singularities$\}$. On a small chart which contains a
singularity, coordinate $z$ can be chosen in such way that
$\psi=z^kdz^2$, where $k$ is the order of the singularity ($k = 0$
corresponds to a regular point, $k=-1$ corresponds to a pole and
$k > 0$ corresponds to a true zero). It is easy to check that in a
neighborhood of a singularity of $\psi$, the metric has a cone
type singularity with the cone angle $(k+2) \pi$.

A meromorphic quadratic differential $\psi$ also defines on $M$ a
pair of transversal foliations. Let $z$ be a canonical coordinate
for $\psi$ i.e. let locally outside of singularity $\psi=(\D
z)^2$. The {\it horizontal} (respectively {\it vertical})
foliation defined by $\psi$ is the foliation $y=\textrm{const}$
(resp $x=\textrm{const}$) locally where $z=x+iy$. The vertical
(horizontal) foliation defined by $\psi$ is oriented if and only
if the quadratic differential $\psi$ is the square of an Abelian
differential.

Similarly, an Abelian differential defines a {\it translation
structure} on the Riemann surface: now the holonomy representation
of the corresponding flat metric in the linear group is trivial.
In particular, all cone angles at the singularities are integer
multiples of $2\pi$.

\begin{figure}[htbp]
\begin{center}
\psfrag{1}{$1$}  \psfrag{4}{$4$}
\psfrag{2}{$2$}  \psfrag{5}{$5$}
\psfrag{3}{$3$}
\includegraphics[width=8cm]{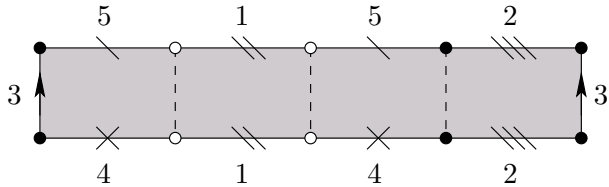}
\end{center}
\caption{
\label{fig:one:cylinder}
Consider the
following identifications of the boundary of this polygon:
identify the corresponding pairs of segments with the numbers
$1,2$ and $3$ by translations; identify the corresponding pairs of
segments with the numbers $4$ and $5$ by centrally symmetries. We
get a half-translation surface. In this case it has genus $g=2$
and the meromorphic differential induced from the quadratic
differential $(\D z)^2$ in the plane has two zeros of order $2$
on this surface. The singularities are conical point in the flat
metric with a cone of angle $4 \pi$. Note that since the
horizontal foliation is not oriented the corresponding quadratic
differential is not the square of an Abelian differential, though
all cone angles of the singularities in the flat metric are
integer multiples of $2\pi$.
}
\end{figure}

Reciprocally, a half-translation structure on a Riemann surface
$M$ and a choice of a distinguished ``vertical'' direction defines
a complex structure and a meromorphic quadratic differential
$\psi$ on $M$. In many cases it is very convenient to present a
quadratic differential with some specific properties by
an appropriate flat surface. Consider, for
example, a polygon in the complex plane $\C$ with the following
property of the boundary: the sides of the polygon are distributed
into pairs, where the sides in each pair are parallel and have
equal length. Identifying the corresponding sides of the boundary
by translations and central symmetries we obtain a Riemann surface
with a natural half-translation flat structure. The quadratic
differential $dz^2$ on $\mathbb{C}$ gives a quadratic differential
on this surface with punctures. The punctures correspond to
vertices of the polygon; they produce the cone type singularities
on the surface. It is easy to see that the complex structure, and
the quadratic differential extends to these points, and that a
singular point of the flat metric with a cone angle $(k+2) \pi$
produces a singularity of order $k$ (a pole if $k=-1$) of the
quadratic differential. (See also Figure~\ref{fig:one:cylinder}
which illustrates this construction.)

\begin{Remark}
The area of $M$ with respect to the flat metric defined by a
meromorphic quadratic differential is finite if and only if $\psi$
does not possess poles of order greater than $1$. This explains
why we consider quadratic differential with simple poles only.
\end{Remark}

\subsection{Non-connectedness of the Strata}

Now, we prove that strata with a hyperelliptic connected component
are not connected in general. We present a particular geometric
property of surfaces, that belong to a hyperelliptic component,
and then we construct appropriate flat surfaces which do not
verify this property.

Let us consider a flat surface $(M,\psi)$ in a hyperelliptic
component of the stratum $\Q(4k+2,4(g-k)-6), \text{ where } k \geq
0, \ g \geq 2, g-k \geq 2$. By definition, there exist a ramified
double covering $\pi : M \rightarrow \CP$ such that $\psi =
\pi^\ast \psi_0$ where $\psi_0$ is a meromorphic quadratic
differential on $\CP$. We consider the canonical atlas on $M$ such
that locally $\psi = \D z^2$ in a neighborhood of a regular point
and $\psi = z^k \D z^2$ in a neighborhood of a singularity. In
this atlas, the hyperelliptic involution $\tau : M \rightarrow M$
is affine. Moreover, $\tau$ is an isometry of the flat metric
defined by $\psi$. The two zeros of $\psi$ are fixed points of
$\tau$.

Suppose that there exist a geodesic saddle connection joining the
two zeros --- a geodesic segment in the flat metric defined by
$\psi$ having no singularities in its interior, and having the two
zeros of $\psi$ as the endpoints. Then this saddle connection have
an image by $\tau$ of the same length and it is also a saddle
connection.

So if there is a horizontal saddle connection (which is
a singular leaf of the horizontal foliation joining  the two
singularities) then there
exists another one of the same length going in the same direction.

In particular, if a quadratic differential $\psi$ having two zeros
of orders $4k+2$ and $4(g-k)-6$ defines a horizontal foliation
such that all horizontal saddle connections between the two zeros
have different lengths, then $(M,\psi)$ does not belong to the
hyperelliptic component of $\Q(4k+2,4(g-k)-6)$. Since the
hyperelliptic component of this stratum is nonempty, it would
imply that the stratum is not connected.

The similar argument can be applied to the strata
$\Q(2(g-k)-3,2(g-k)-3,2k+1,2k+1)$ and
$\Q(2(g-k)-3,2(g-k)-3,4k+2)$.

We use this idea to construct appropriate flat surfaces such that
the corresponding quadratic differentials do not belong to
hyperelliptic components. We use the following two Lemmas which
are particular cases of the corresponding Lemmas of Eskin, Masur
and Zorich (\cite{Eskin:Masur:Zorich},
\cite{Eskin:Masur:Zorich:2}); see also analogous Proposition
$4.7(b)$ in the paper~\cite{Hubbard:Masur} of Hubbard and Masur.

\begin{Lemma}
\label{lm:separation}
Consider a surface in $\Q(k_1,\dots,k_n)$. Choose $l_1,l_2 \in
\{-1,1,2,3,4,\dots\}$, as follows
\begin{itemize}
\item if $k_1$ is odd, $l_1+l_2 = k_1$, $l_i$ any.
\item if $k_1$ is even, $l_1+l_2 = k_1$, $l_i$ even.
\end{itemize}
For any $\psi_0\in\Q(k_1,k_2,\dots,k_n)$ and for any sufficiently
small $\varepsilon > 0$ (depending on $\psi_0$) it is possible to
construct a deformation $\psi\in\Q(l_1,l_2,k_2,\dots,k_n)$ of
$\psi_0$ such that the corresponding flat metric has a
horizontal saddle connection of length $\varepsilon$ joining the
singularities $P_1$ and $P_2$ of orders $l_1$ and $l_2$.

The deformation can be chosen to be local: the flat metric does
not change outside of a small neighborhood of the zero of
multiplicity $k_1$.
\end{Lemma}
\begin{Lemma}
\label{lm:separation:in:three}
Consider a surface in $\Q(k_1,\dots,k_n)$. Let $k_1$ be odd and
let $k_1=l_1+l_2+l_3$, where $l_i\in \{-1,1,2,3,4,\dots\}$ are
also odd. For any $\psi_0\in\Q(k_1,k_2,\dots,k_n)$ and for any
sufficiently small $\varepsilon > 0$ (depending on $\psi_0$) it is
possible to construct a deformation
$\psi\in\Q(l_1,l_2,l_3,k_2,\dots,k_n)$ of $\psi_0$ such that the
corresponding flat metric has two horizontal saddle
connection of length $\varepsilon$ joining the singularities
$P_1$, $P_2$ and $P_2, P_3$ of orders $l_1, l_2, l_3$
correspondingly.

The deformation can be chosen to be local: the flat metric does
not change outside of a small neighborhood of the zero of
multiplicity $k_1$.
\end{Lemma}
\begin{proof}[Proof of Lemma~\ref{lm:separation}]
Let $P_0\in M$ be the singularity of $\psi_0$ of order $k_1$.
Consider a small metric disk $D$ with the center in $P_0$ and of
radius $R$ in $M$ in the flat metric on $M$ defined by the
quadratic differential $\psi_0$. We choose $R$ to be small enough,
so that $D$ does not contain other singularities of $\psi_0$.

In fact, $D$ is glued from $k_1+2$ Euclidean half-disks of radius
$R$ where the corresponding radii of the half-disks are pairwise
identified, see the left sides of Figures~\ref{break:4:22}
and~\ref{break:3:12}. Note that the pictures are schematic: the
angle of every sector is actually equal to $\pi$.

To make a local deformation of the Euclidean metric inside $D$ we
reglue the radii borders of Euclidean half-disks in a different
way. Figure~\ref{break:4:22} illustrates how to break a zero of
even degree $k_1$ when even number $k_1+2$ of half-disks are
adjacent to the vertex into two zeros of even degrees $l_1$ and
$l_2$. Now there are $l_1+2$ and $l_2+2$ half-disks adjacent to
the corresponding vertices. Note that the angles of all sectors
are, actually, again equal to $\pi$. Figure~\ref{break:3:12}
illustrates how to break a zero of odd degree $k_1$ into a zero of
even degree $l_1$ and a zero of odd degree $l_2$.

\begin{figure}[htbp]
\begin{center}
\psfrag{6p}{\scriptsize $6\pi$}
\psfrag{e}{\scriptsize $\varepsilon$}
\psfrag{d}{\scriptsize $2\delta$}
\psfrag{e+}{\scriptsize $\varepsilon+\delta$}
\psfrag{e-}{\scriptsize $\varepsilon-\delta$}

\includegraphics[width=13cm]{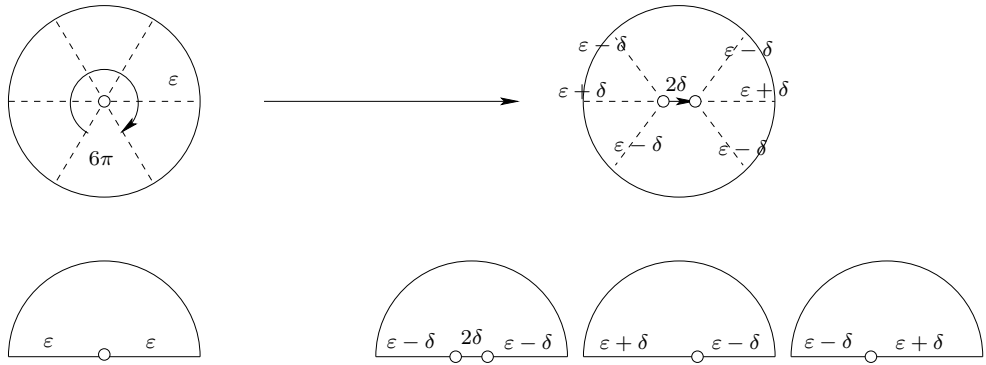}
\end{center}
\caption{
\label{break:4:22}
Breaking up a zero of order $4$ into two zeros of orders $2$.
Note that the surgery is local: we do not change the
flat metric outside of the neighborhood of the zero.}
\end{figure}

\begin{figure}[htbp]
\begin{center}
\psfrag{6p}{\scriptsize $5\pi$}
\psfrag{e}{\scriptsize $\varepsilon$}
\psfrag{d1}{\scriptsize $\delta$}
\psfrag{e+}{\scriptsize $\varepsilon+\delta$}
\psfrag{e-}{\scriptsize $\varepsilon-\delta$}

\includegraphics[width=13cm]{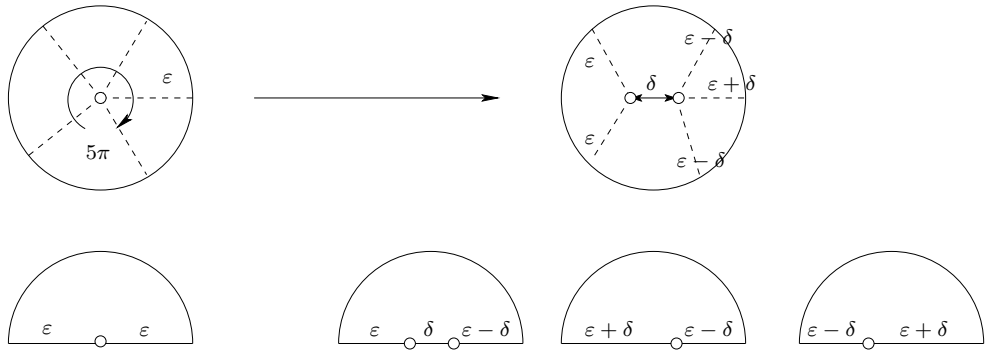}
\end{center}
\caption{
\label{break:3:12}
Breaking up a zero of order $3$ into two zeros of orders $1$ and
$2$ correspondingly. Note that the surgery is local: we do not change the
flat metric outside of the neighborhood of the zero.}
\end{figure}

\end{proof}

\begin{figure}[htbp]
\begin{center}
\psfrag{6p}{\scriptsize $5\pi$}
\psfrag{e}{\scriptsize $\varepsilon$}
\psfrag{d1}{\scriptsize $\delta$}
\psfrag{d2}{\scriptsize $\varepsilon-2\delta$}
\psfrag{e+}{\scriptsize $\varepsilon+\delta$}
\psfrag{e-}{\scriptsize $\varepsilon-\delta$}

\includegraphics[width=13cm]{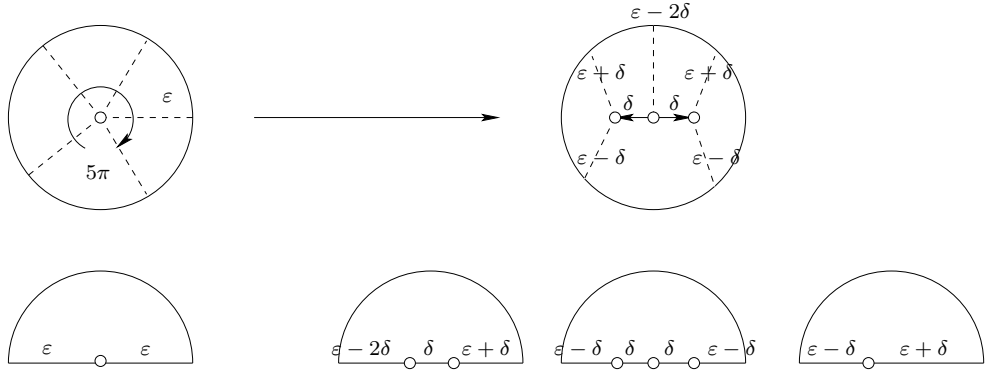}
\end{center}
\caption{
\label{fig:separation:in:three}
Breaking up a zero of order $3$ into three zeros of orders $1$.
Note that the surgery is local: we do not change the flat metric
outside of the neighborhood of the zero.}
\end{figure}

\begin{proof}[Proof of Lemma~\ref{lm:separation:in:three}]
The proof of Lemma~\ref{lm:separation:in:three} is completely
analogous to the proof of Lemma~\ref{lm:separation}. It is
illustrated at figure~\ref{fig:separation:in:three}.
\end{proof}

\begin{Theorem}
\label{th:non:connected}
The strata which possess a hyperelliptic component 
are non connected except for the following cases 
corresponding respectively to genus $1$ and $2$
\begin{displaymath}
\begin{array}{ccc}
\begin{cases}
\Q(-1,-1,2) \\
\Q(-1,-1,1,1)
\end{cases}  &
and &
\begin{cases}
\Q(2,2) \\
\Q(1,1,2) \\
\Q(1,1,1,1)
\end{cases}
\end{array}
\end{displaymath}
when the whole stratum coincides with its hyperelliptic connected
component.
\end{Theorem}
\begin{proof}[Proof of Theorem~\ref{th:non:connected}]
We decompose the proof into three Lemmas. In the first, we prove that the strata
of the list above are connected. In the second, we prove Theorem~\ref{th:non:connected} for
strata with $2$ and $3$ singularities. In the third we prove the Theorem for
strata with fours singularities.
\begin{Lemma}
\label{lm:p:1}
The following strata are connected
$$\Q(-1,-1,2), \Q(-1,-1,1,1) \textrm{ in genus $1$}$$
$$\Q(2,2), \Q(1,1,2), \Q(1,1,1,1) \textrm{ in genus $2$}$$
they coincide with the corresponding hyperelliptic component.
\end{Lemma}
\begin{proof}[Proof of Lemma~\ref{lm:p:1}]
We use an idea drawn from the book of Farkas-Kra~\cite{Farkas:Kra}.
See also~\cite{Lanneau:02:classification} for another proof, using the combinatorics of 
``generalized permutations''.

We prove this Lemma for $\Q(2,2)$. For the other strata,
the demonstration is similar. First of all, take a point
$[M,\psi] \in \Q(2,2)$. In genus $2$, all curves are hyperelliptic
so we can suppose that a representative of $M$ is given by an
algebraic curve$$w^2 = \prod_{i=1}^6(z-z_i) \qquad (w,z) \in \C^2$$
where $z_i \not = z_j$ for $i \not = j$ and $z_i \in \C$. There
is a (ramified) double covering $\pi : M \rightarrow \CP$ defined
by $\pi(w,z)=z$. Let $\tau$ be the hyperelliptic involution
$\tau(w,z)=(-w,z)$. With these notations, we can see
that $\omega_1=\cfrac{\D z}{w}$ and $\omega_2 =z \cdot \cfrac{\D z}{w}$
form a basis of holomorphic Abelian differentials on $M$. It easy to see
that $\omega_1^2,\omega_2^2$ and $\omega_1 \cdot \omega_2$ are linearly
independent in the space of holomorphic quadratic
differentials which has dimension $3g-3=3$. So the family $\{
\omega_1^2,\omega_2^2,\omega_1 \cdot \omega_2 \}$ forms a basis
in the space of holomorphic quadratic differentials on $M$.
By direct computation, we can see that $\tau^\ast \omega_i=- \omega_i$
and consequently, any {\it holomorphic} quadratic differential
is invariant under $\tau$.

Let $(M,\psi)\in \Q(2,2)$ and let
$P_1$ and $P_2$ be the two zeros of $\psi$. We have
$$
0=\psi(P_1)=\tau^\ast\psi(P_1)=\psi(\tau(P_1))
$$
hence $\tau$ either fixes the two zeros or interchanges them.

Suppose that the second case occurs. Let $z_0=\pi(P_1)=\pi(P_2)$.
Since by assumption $\tau(P_1)\neq P_2$, we see that $z_0$ is a
noncritical value of the covering, and hence $z_0\neq z_i$ for any
$i=1, \dots, 6$. We can construct a quadratic differential
$\psi_0$ on $\CP$ with six simple poles at the images $z_i \in
\CP$ of the ramifications points, and a zero of order $2$ at
$z_0$. Then the covering $\pi : M \rightarrow \CP$ is exactly the
covering of Construction~\ref{constr:canonical:2:covering} so
$\pi^\ast\psi_0=\omega^2$ possesses two zeros of order two at the
same points as the two zeros of $\psi$. Hence $\psi=const
\cdot \omega^2$ which contradicts the assumption that $\psi$ is
not a global square of an Abelian differential.

Thus the two zeros $P_1$ and $P_2$ are the fixed points of the
hyperelliptic involution. Let $z_0=\pi(P_1),z^{'}_0=\pi(P_2)$. We
can construct a quadratic differential $\psi_0$ on $\CP$ with four
poles at the images of ramifications points, and with two fake
zeros at $z_0$ and $z^{'}_0$ correspondingly. As above,
$\psi=const \cdot \pi^\ast\psi_0$. But by construction
$\pi^\ast\psi_0$ belongs to the hyperelliptic connected component
(see construction~\ref{constr:general}). Thus $p \in
\Q^{hyp}(2,2)$ and hence $\Q(2,2)=\Q^{hyp}(2,2)$ is connected.

This completes the proof of Lemma~\ref{lm:p:1} in the case of the
stratum $\Q(2,2)$. The proof for other cases is completely
analogous.
\end{proof}
Now we prove Theorem~\ref{th:non:connected} for the general strata with $2$ and $3$
singularities.
\begin{Lemma}
\label{lm:p:2}
The strata
$$
\Q(4(g-k)-6,4k+2) \textrm{ with } k \geq 0, g \geq 3 \textrm{ and }
g - k \geq 2
$$
and
$$\Q(2(g-k)-3,2(g-k)-3,4k+2)$$
$$\textrm{ with } g \geq 3, k\geq 0, g-k\geq 1
\textrm{ or } g=2 \textrm{ and } k=1$$
are non-connected.
\end{Lemma}
\begin{proof}[Proof of Lemma~\ref{lm:p:2}]
By the Theorem of Masur and Smillie~\cite{Masur:Smillie} (see also
the statement of the theorem at the beginning of
section~\ref{s:covering:construction}) the stratum $\Q(4g-4)$ is
non-empty when genus $g \geq 3$. Consider a flat surface
corresponding to some $[M,\psi_0]\in \Q(4g-4)$
such that the horizontal measured foliation is uniquely ergodic,
in particular minimal. (The genericity of this property is proved
in~\cite{Masur:82} and in~\cite{Veech:82}.) Consider a small disk
around the zero of the quadratic differential. Applying Lemma
above, we can break the zero of order $4g-4$ into two zeros of
orders $4(g-k)-6$ and $4k+2$ correspondingly with a very short
horizontal saddle connection $\gamma$ joining them. By
construction there are no other short horizontal saddle
connections. So the quadratic
differential $\psi\in\Q(4(g-k)-6,4k+2)$ thus constructed does not
belong to hyperelliptic connected component of
$\Q(4(g-k)-6,4k+2)$ (see the beginning of this section).
This proves the first statement.


Consider now the stratum $\Q(2(g-k)-3,2(g+k)-1)$ with $k \geq 0$,
$g-k\geq 1$. By the Theorem of Masur and
Smillie~\cite{Masur:Smillie} this stratum is non-empty if $g \geq
3$ or if $g=2$ and $k=1$.

Take a flat surface corresponding to a quadratic differential from
this stratum such that the corresponding horizontal foliation is
minimal (see above). Then if we break the zero of order $2(g+k)-1$
into two zeros of orders $2(g-k)-3$ and $4k+2$ joined by a very
short horizontal saddle connection we get a surface which belongs
to the non-hyperelliptic component of the stratum
$\Q(2(g-k)-3,2(g-k)-3,4k+2)$.

The only strata to which we cannot apply this method are the
stratum with two poles and a simple zero $\Q(-1,-1,2)$ in genus 
$1$ and the stratum with two simples zeros and a double zero 
$\Q(1,1,2)$ in genus $2$. In these two cases the initial 
(non-perturbed) quadratic differentials would belong to 
$\Q(-1,1)$ and $\Q(1,3)$ which are empty.

Lemma~\ref{lm:p:2} is proved.
\end{proof}
Let us prove proposition for strata with four singularities.
\begin{Lemma}
\label{lm:p:3}
The strata
$$\Q(2(g-k)-3,2(g-k)-3,2k+1,2k+1)$$
$$\textrm{ with } g \geq 3, k\geq -1,g-k\geq 2
\textrm{ or } g=2 \textrm{ and } k=-1$$ are non-connected.
\end{Lemma}
\begin{proof}[Proof of Lemma~\ref{lm:p:3}]
Consider the stratum $\Q(2(g-k)-3,2(g+k)-1)$ with $g\geq 3$,
$k \geq -1$, $g-k\geq 2$. By the Theorem of Masur and
Smillie~\cite{Masur:Smillie} this stratum is non-empty.

Take a flat surface corresponding to a quadratic differential from
this stratum such that the corresponding horizontal foliation is
minimal (see above). By Lemma~\ref{lm:separation:in:three}
we can break the zero of order $2(g+k)-1$
into three zeros of orders $2(g-k)-3, 2k+1, 2k+1$,
joined consecutively by  very
short horizontal saddle connections. We use a local perturbation
of the flat metric, which not rather create a short horizontal saddle
connection going to the unperturbed zero of order $2(g-k)-3$.
We get a flat surface corresponding to a quadratic differential
from the stratum $\Q(2(g-k)-3,2(g-k)-3,2k+1,2k+1)$.

The hyperelliptic involution interchanges
the corresponding zeros of any
$\psi'\in\Q^{hyp}(2(g-k)-3,2(g-k)-3,2k+1,2k+1)$.
By construction the quadratic differential constructed
above is asymmetric, and hence it belongs
to the non-hyperelliptic component of the stratum
$\Q(2(g-k)-3,2(g-k)-3,4k+2)$.
Lemma~\ref{lm:p:3} is proved.
\end{proof}
Now, Theorem~\ref{th:non:connected} follows from
Lemmas~\ref{lm:p:2} and~\ref{lm:p:3}.
\end{proof}

\section{Announcement of the Classification Theorem}
\label{general:theorem}

In~\cite{Kontsevich:Zorich} Kontsevich and Zorich have shown that the   
connected components of the moduli spaces of Abelian differentials
are classified exactly by two invariants: the hyperellipticity and the 
parity of the spin structure. 

The story for the moduli space of quadratic differential is more complicated. 
For example, the strata $\Q(12)$ and $\Q(-1,9)$ do not have any hyperelliptic 
components. In addition, in paper~\cite{Lanneau:02:spin}, we show that all 
quadratic differentials in a fixed strata must have the same parity of the spin structure. 
In the previous case, the spin structure is $0$ for any quadratic differential in each
of these two strata. However, it was proved by A.~Zorich by a direct
computation of the corresponding {\it extended Rauzy classes} that
each of these two strata has exactly two distinct connected components.

In~\cite{Lanneau:02:classification} we give the following general
description of all connected components of any stratum of the moduli
spaces $\Q(k_1,\dots,k_n)$:
     
\begin{Theorem}
Let us fix $g \geq 5$. Then all strata of the moduli space $\Q_g$ 
listed in Theorem~\ref{th:non:connected} have exactly two connected components: 
one is hyperelliptic --- the other not.

All other strata of the moduli space of meromorphic quadratic differentials 
$Q_g$ are non-empty and connected.
\end{Theorem}

And for small genera we have:

\begin{Theorem}
Let us fix $g \leq 4$. The components of the strata of the moduli
space $\Q_g$ fall in the following description
\begin{itemize}

\item In genera $0$ and $1$, any stratum is connected.

\item In genus $2$ there are two hyperelliptic non-connected
strata. All other strata of $\Q_2$ are connected.

\item In genera $3$ and $4$, any hyperelliptic stratum possesses
two connected components: one is hyperelliptic --- the other not.
All other strata, with $4$ exceptional cases, are connected.

\item The $4$ above particular cases are
$$
\Q_{g=3}(-1,9), \ \Q_{g=3}(-1,3,6), \ \Q_{g=3}(-1,3,3,3), \ \Q_{g=4}(12)
$$
and these strata have two connected components.

\end{itemize}
\end{Theorem}

\section{Acknowledgments}
I would like to thank A.~Zorich for his proposition to work with
such beautiful objects like quadratic differentials and moduli
spaces. I am grateful to A.~Zorich, P.~Hubert and the Referee for their
remarks concerning this paper.

\end{document}